\documentclass[12pt]{amsart}
\usepackage{amsmath}
\usepackage{amssymb}
\usepackage{amsthm}
\usepackage{amscd}

\usepackage{color}
\usepackage{hyperref}

\usepackage{tikz-cd}
\usetikzlibrary{matrix} % matriz no tikz
\usetikzlibrary{cd} % diagrama comutativa

\usepackage[latin1]{inputenc}
\usepackage{mathpazo}
\usepackage[scaled=.95]{helvet}
\usepackage{courier}
\usepackage[all]{xy}

\swapnumbers
\headsep=1cm \numberwithin{equation}{section}
\usepackage[colorinlistoftodos]{todonotes}

\newtheorem{theorem}{Theorem}[section]

\newtheorem{lemma}[theorem]{Lemma}

\newtheorem{proposition}[theorem]{Proposition}

\newtheorem{corollary}[theorem]{Corollary}

\theoremstyle{definition}

\newtheorem{definition}[theorem]{Definition}

\newtheorem{example}[theorem]{Example}

\newtheorem{remark}[theorem]{Remark}

\numberwithin{equation}{section}

\newtheorem{remark/questions}[theorem]{Remark and Questions}
\newtheorem{fact}[theorem]{Fact}

\newtheorem*{fund}{Funding}
\newtheorem*{agra}{Acknowlegment}

 \long\def\alert#1{\smallskip\line{\hskip\parindent\vrule
\vbox{\advance\hsize-2\parindent\hrule\smallskip\parindent.4\parindent
  \narrower\noindent#1\smallskip\hrule}\vrule\hfill}\smallskip}

\newtheorem{newclaim}[theorem]{}

\def\dim{\mathop{\rm dim}}

\numberwithin{equation}{section}

\begin{document}

\dedicatory{}

\title[]{On $\operatorname{Ext}$-finite modules, quasi-injective dimension and width of modules}
%A sufficent condition for the positivity of the Huneke-Wiegand conjecture
\author{V. H. Jorge-P\'erez}

\address{Universidade de S{\~a}o Paulo -
ICMC, Caixa Postal 668, 13560-970, S{\~a}o Carlos-SP, Brazil}
\email{vhjperez@icmc.usp.br}

\author{Paulo Martins}

\address{Universidade de S{\~a}o Paulo -
ICMC, Caixa Postal 668, 13560-970, S{\~a}o Carlos-SP, Brazil}
\email{paulomartinsmtm@gmail.com}

\thanks{Corresponding author: Paulo Martins}

\begin{abstract}
Let $(R,\mathfrak{m},k)$ be a commutative Noetherian local ring. It is well-known that if $M$ is a finitely generated $R$-module of finite quasi-injective dimension, then $\operatorname{qid}_RM=\operatorname{depth}R$. In this paper, we demonstrate that under the weaker condition that $M$ is $\operatorname{Ext}$-finite and has finite quasi-injective dimension, the equality $\operatorname{qid}_RM=0$ holds if and only if $\operatorname{Ext}_R^{i>0}(R/(\boldsymbol{x}),M)=0$, where $\boldsymbol{x}\in\mathfrak{m}$ is a maximal $R$-sequence and if $\operatorname{qid}_RM\neq0$, we show then that $\operatorname{qid}_RM=\sup\lbrace i:\operatorname{Ext}_R^i(R/(\boldsymbol{x}),M) \neq0\rbrace$. Also, we show that if $R$ is a Cohen-Macaulay local ring and $M$ is an $\operatorname{Ext}$-finite $R$-module of finite quasi-injective dimension, then $\operatorname{depth}R=\operatorname{qid}_R M+\inf\lbrace i:\operatorname{Tor}_i^R(k,M)\neq0\rbrace$, provided that $\inf\lbrace i:\operatorname{Tor}_i^R(k,M)\neq0\rbrace<\infty$. 
\end{abstract}

\maketitle

\section{Introduction}
Throughout this paper, $R$ will denote a commutative Noetherian ring. The quasi-injective dimension (qid) is a refinement of the classical notion of the injective dimension of a module, in the sense that there is always an inequality $\operatorname{qid}_RM \leq  \operatorname{id}_R M$ for all $R$-modules $M$. It was introduced by Gheibi \cite{Gheibi} that recovered several well-known results about injective and Gorenstein-injective dimensions in the context of quasi-injective dimension. For a finitely generated $R$-module $M$ of finite quasi-injective dimension over a local ring, Gheibi established a version of the Bass formula, proving that $\operatorname{qid}_R M = \operatorname{depth} R$. Also, in a recent paper, Tri \cite{Tri} established a version of Chouinard's formula for quasi-injective dimension. However, a formula for the quasi-injective dimension has not been proved when $R$ is local and $M$ is not finitely generated.

Enochs and Jenda \cite{ENOCHS}, Sahandi and Sharif \cite{SHARIF}, as well as Sharif and Yassemi \cite{YASSEMI}, among others, established partial generalizations of the Bass formula under the weaker assumption that the $R$-module $M$ is \(\operatorname{Ext}\)-finite. That is, they assumed that for every finitely generated \( R \)-module \( N \), the $R$-module \( \operatorname{Ext}_R^i(N,M) \) is finitely generated for all \( i \geq 1 \). Their results focused on the cases where \( M \) has finite Gorenstein-injective dimension or finite injective dimension and established conditions in which versions of the following formula are satisfied:
\begin{align}\label{equation}
   \operatorname{id}_R M + \operatorname{width}_R M  \; = \operatorname{depth} R .
\end{align}
Recall that the width of a module \( M \) over a local ring \( R \) with residue field \( k \) is defined as 
\(
\operatorname{width}_R M = \inf \lbrace i : \operatorname{Tor}_i^R(k, M) \neq 0 \rbrace\).

Motivated by these articles, the objective of this paper is to provide a homological characterization of quasi-injective dimension and to establish a version of equality (\ref{equation}), both in the context of \(\operatorname{Ext}\)-finite modules of finite quasi-injective dimension. Let us now describe the contents of the paper. In Section 2, we begin by introducing the necessary notations, definitions, and foundational results that will be utilized throughout the paper. Section 3 is devoted to proving results concerning \(\operatorname{Ext}\)-finite modules of finite quasi-injective dimension. Under these conditions, we derive a homological characterization of the quasi-injective dimension. The main results of this section are stated below:

\begin{theorem}\label{prop1}{\rm (See Theorem \ref{prop1})}
Let $(R,\mathfrak{m},k)$ be a local ring and let $M$ be a non-zero $\operatorname{Ext}$-finite $R$-module of finite quasi-injective dimension. The following conditions are equivalent.
\begin{enumerate}
    \item $\operatorname{qid}_R M \leq n$. 
    \item  $\operatorname{Ext}_R^{i>n}(N,M)=0$ for all $R$-modules $N$ such that $\operatorname{Ext}^{i \gg 0}_R(N,M)=0.$
    \item $\operatorname{Ext}^{i>n}_R(R/J,M)=0$ for all ideals $J$ of $R$ such that $\operatorname{Ext}_R^{i\gg0} (R/J,M)=0.$
    \item $\operatorname{Ext}^{i>n}_R(R/(\boldsymbol{x}),M) = 0$, if $\boldsymbol{x} \in \mathfrak{m}$ is a maximal $R$-sequence.   
\end{enumerate}
\end{theorem}

\begin{corollary}{\rm (See Corollary \ref{crl:bom})}
    Let \((R, \mathfrak{m},k)\) be a local ring and let \(M\) be a non-zero \(\operatorname{Ext}\)-finite \(R\)-module of finite quasi-injective dimension. Then, for any maximal \(R\)-sequence \(\boldsymbol{x} \in \mathfrak{m}\), the following holds.
\begin{enumerate}
    \item $\operatorname{qid}_R M =0$ if and only if $\operatorname{Ext}_R^{i>0}(R/(\boldsymbol{x}),M)=0.$
    \item If $\operatorname{qid}_R M \neq 0$, then         $\operatorname{qid}_R M = \sup \lbrace i:\operatorname{Ext}_R^i(R/(\boldsymbol{x}),M ) \neq 0\rbrace$. 
\end{enumerate}
\end{corollary}
In Section 4, we present our main theorem, which establishes a version of equality (\ref{equation}) for quasi-injective dimension. Specifically, we prove the following result:
\begin{theorem}\label{theorem1} \rm{(See Theorem \ref{theorem1})}
Let $(R,\mathfrak{m},k)$ be a Cohen-Macaulay local ring and let $M$ be a non-zero $\operatorname{Ext}$-finite $R$-module with $\operatorname{width}_R M<\infty$. If $M$ is of finite quasi-injective dimension, then $$\operatorname{qid}_R M + \operatorname{width}_R M = \operatorname{depth} R.$$ 
\end{theorem}

\section{Setup and background}
In this section, we introduce fundamental definitions and facts, such as quasi-injective dimension, cosyzygyes and $\operatorname{Ext}$-finite modules. These
concepts will be essential throughout the rest of the paper.
\begin{newclaim}    
For a complex $$X: \cdots \stackrel{\partial_{i+2}}{\longrightarrow} X_{i+1} \stackrel{\partial_{i+1}}{\longrightarrow  } X_i \stackrel{\partial_{i}}{\longrightarrow
      } X_{i-1}   \longrightarrow  \cdots $$ of $R$-modules, we set for each integer $i$, $\operatorname{Z}_i(X)=\ker \partial_i$ and $\operatorname{B}_i(X)= \operatorname{Im} \partial_{i+1}$ and $\operatorname{H}_i(X)=\operatorname{Z}_i(X)/\operatorname{B}_i(X)$.  Moreover, we set: \begin{align*}
&\left\{ 
    \begin{aligned}
        \sup X &= \sup \{ i \in \mathbb{Z} : X_i \neq 0 \},\\
        \inf X &= \inf \{ i \in \mathbb{Z} : X_i \neq 0 \},
    \end{aligned}
\right. \quad 
\left\{
    \begin{aligned}
        \text{hsup } X &= \sup \{ i \in \mathbb{Z} : \operatorname{H}_i(X) \neq 0 \},\\
        \text{hinf } X &= \inf \{ i \in \mathbb{Z} : \operatorname{H}_i(X) \neq 0 \}.
    \end{aligned}
\right.
\end{align*}
The \textit{length} of $X$ is defined to be  $\operatorname{lenght} X= \sup X - \inf X$. We say that $X$ is \textit{bounded}, if $\operatorname{length} X < \infty$. We say that $X$ is \textit{bounded above} if $\sup X < \infty$.
\end{newclaim}
The definition of quasi-injective dimension was introduced by Gheibi \cite{Gheibi} as a generalization of the classical notion of injective dimension.
\begin{definition}
    Let $M$ be an $R$-module. A \textit{quasi-injective resolution} of $M$ is a bounded above complex $I$ of injective $R$-modules such that for all $i \leq \sup I$ there exist non-negative integers $b_i$, not all zero, such that $H_i(I) \cong M^{\oplus b_i}$. The \textit{quasi-injective dimension} of $M$ is defined by 
  \begin{align*}
            \operatorname{qid}_R M = \inf \lbrace \operatorname{hinf} I - \operatorname{inf} I :  I \text{ is a bounded quasi-injective resolution of } M \},
        \end{align*}
        if  $M\not=0$, and $\operatorname{qid}_R M=-\infty$ if $M=0$.
    One has $\operatorname{qid}_R M=\infty$ if and only if $M$ does not admit a bounded quasi-injective resolution. 
\end{definition}

We remark that quasi-injective dimension is finer than injective dimension, in the sense that there is always an inequality $\operatorname{qid}_R M \leq \operatorname{id}_R M$. Also, Gheibi proved that equality holds if $\operatorname{id}_R M < \infty$ and $M$ is a finitely generated module over a local ring (see \cite[Corollary 3.3]{Gheibi}).

\begin{fact}\label{newfact}
(\cite[Proposition 3.4(2)]{Gheibi}) Let $M,N$ be $R$-modules and assume that $\operatorname{qid}_R M < \infty$. If $\operatorname{Ext}_R^i(N,M)=0$ for $i \gg 0$, then $\operatorname{Ext}_R^i(N,M)=0$ for all $i> \operatorname{qid}_R M.$
\end{fact}

\begin{newclaim}
Let $M$ be an $R$-module. Consider the augmented minimal injective resolution of $M$:
\begin{align*}
0 \rightarrow M \xrightarrow{d^{-1}} E(M)  \xrightarrow{d^0} E^1(M) \xrightarrow{d^1} E^2(M) \rightarrow \cdots \rightarrow E^{n-1}(M) \xrightarrow{d^{n-1}} E^n(M) \rightarrow \cdots 
\end{align*}
Recall that the $n$-th \textit{cosyzygy} of $M$ is defined by $\Omega_{-n}^RM:= \operatorname{Im}(d^{n-1})$ for all $n \geq 0$.    
\end{newclaim}

Enochs and Jenda \cite{ENOCHS} introduced the concept of \(\operatorname{Ext}\)-finite modules. The $R$-module $M$ is called $\operatorname{Ext}$-{\it finite} if the $R$-module $\operatorname{Ext}_R^i(N,M)$ is finitely gene\-rated for every finitely generated $R$-module $N$ and for all $i \geq 1$. Clearly, every finitely generated $R$-module is $\operatorname{Ext}$-finite. Moreover, it is straightforward to verify that every cosyzygy of an $\operatorname{Ext}$-finite $R$-module remains $\operatorname{Ext}$-finite. 

\begin{lemma}{\rm (\cite[Lemma 2.2]{ENOCHS})}\label{lemma:enochs} Let $(R,\mathfrak{m},k)$ be a local ring and let $M$ be an $\operatorname{Ext}$-finite $R$-module. Then 
        $\operatorname{id}_R M = 
\sup \lbrace i:\operatorname{Ext}_R^i(k,M) \neq 0 \rbrace$.
\end{lemma}

\begin{theorem} {\rm (\cite[Theorem 2.6]{ENOCHS})}\label{finite} Let $R$ be a $d$-dimensional local ring and let $M$ be a non-injective $\operatorname{Ext}$-finite $R$-module. Then $\operatorname{width}_R M \leq d$.
\end{theorem}
\begin{proposition}{\rm (\cite[Proposition 2.7]{ENOCHS})}\label{vanishing} Let $(R,\mathfrak{m},k)$ be a local ring and let $M$ be an injective $R$-module. Then $\operatorname{width}_RM < \infty$ if and only if $E(k)$ is a summand of $M$. In this case, $\operatorname{width}_R M = \operatorname{depth} R$.
    
\end{proposition}
The next lemma is a reformulated version of a lemma of Gheibi, replacing the condition of \( M \) being finitely generated with the condition that \( M \) is \( \operatorname{Ext} \)-finite.
\begin{lemma}\label{lemma:similar}
Let $M$ be a non-zero $\operatorname{Ext}$-finite $R$-module of finite quasi-injective dimension and let
$$I: \,\, 0 \rightarrow I_0 \xrightarrow{\partial_0} I_{-1}  \xrightarrow{\partial_{-1}} I_{-2} \rightarrow \cdots  $$
be a bounded quasi-injective resolution of $M$ such that $\operatorname{qid}_R M = \operatorname{hinf} I - \operatorname{inf} I$ and  $\operatorname{sup} I= \operatorname{hsup} I=0$ (such $I$ always exists, by \cite[Remark 2.3(3)]{Gheibi}). Set $Z_i=Z_i(I)$ and $B_i=B_i(I)$, for all $i \in \mathbb{Z}$. Then $B_i$ and $Z_i$ are $\operatorname{Ext}$-finite $R$-modules for all $i \in \mathbb{Z}$.
\end{lemma}
For the proof of the previous lemma, follow the standard argument introduced by Gheibi (see \cite[Lemma 3.1]{Gheibi}), but replace the assumption that \( M \) is a finitely generated \( R \)-module with the weaker condition that \( M \) is an \(\operatorname{Ext}\)-finite \( R \)-module.

\section{Quasi-injective dimension and $\operatorname{Ext}$-finite modules}
In this section, we present a homological characterization of the quasi-injective dimension, under the assumption that the $R$-module \( M \) is \( \operatorname{Ext} \)-finite. The next theorem is comparable to a classical result for injective dimension (see \cite[Proposition 3.1.10]{BH}).

\begin{theorem}\label{prop1}
Let $(R,\mathfrak{m},k)$ be a local ring, and let $M$ be a non-zero $\operatorname{Ext}$-finite $R$-module of finite quasi-injective dimension. The following conditions are equivalent.
\begin{enumerate}
    \item $\operatorname{qid}_R M \leq n$. 
    \item  $\operatorname{Ext}_R^{i>n}(N,M)=0$ for all $R$-modules $N$ such that $\operatorname{Ext}^{i \gg 0}_R(N,M)=0.$
    \item $\operatorname{Ext}^{i>n}_R(R/J,M)=0$ for all ideals $J$ of $R$ such that $\operatorname{Ext}_R^{i\gg0} (R/J,M)=0.$
    \item $\operatorname{Ext}^{i>n}_R(R/(\boldsymbol{x}),M) = 0$, if $\boldsymbol{x} \in \mathfrak{m}$ is a maximal $R$-sequence.   
\end{enumerate}
\end{theorem}
\begin{proof}
(1) $\Rightarrow$ (2) Follows from Fact \ref{newfact}. 

(2) $\Rightarrow$ (3) and (3) $\Rightarrow$ (4) are trivial.

(4) $\Rightarrow (1)$ Consider the notation previously introduced in Lemma \ref{lemma:similar} and let $\boldsymbol{x} \in \mathfrak{m}$ be a maximal $R$-sequence. Set  $s=\operatorname{hinf} I$. We claim that $\operatorname{Ext}_{R}^{i>n}(R/(\boldsymbol{x}),Z_s)=0$. There are exact sequences
\begin{align}\label{sequences}
\begin{cases}
0 \rightarrow Z_j \rightarrow I_j \rightarrow B_{j-1} \rightarrow 0 \\  0 \rightarrow B_j \rightarrow Z_j \rightarrow H_j(I) \rightarrow 0
\end{cases}  
\quad (j \in \mathbb{Z}).
\end{align}
with $H_j(I) \cong M^{\oplus b_j}$ for some $b_j \geq 0$. Consider the exact sequence $0 \rightarrow Z_0 \rightarrow I_0 \rightarrow B_{-1} \rightarrow 0$ and note that $Z_0 \cong \oplus^{b_0} M$ for some positive integer $b_0$. This exact sequence induces the following exact sequence:

\begin{align*}
 \cdots \rightarrow \operatorname{Ext}_R^i(R/(\boldsymbol{x}),&M)^{\oplus b_0} \rightarrow  \operatorname{Ext}_R^{i}(R/(\boldsymbol{x}),I_0)  
 \\ 
& \rightarrow\operatorname{Ext}_{R}^i(R/(\boldsymbol{x}),B_{-1}) \rightarrow \operatorname{Ext}_R^{i+1}\left(R/(\boldsymbol{x}),M\right)^{\oplus b_0} \rightarrow \cdots   .
\end{align*}
Thus, we see that $\operatorname{Ext}_R^{i>n}(R/(\boldsymbol{x}),B_{-1})=0$. Using this vanishing, considering the exact sequence $0 \rightarrow B_{-1} \rightarrow Z_{-1} \rightarrow H_{-1}(I) \rightarrow 0$ and using again the long exact sequence $\operatorname{Ext}_R(R/(\boldsymbol{x}),-)$ we see that $\operatorname{Ext}_R^{i>n}(R/(\boldsymbol{x}),Z_{-1})=0$. 
Therefore, by using the exact sequences (\ref{sequences}) and repeating this argument a finite number of times, we conclude that \(\operatorname{Ext}^{i>n}_R(R/(\boldsymbol{x}), Z_s) = 0\).

Now, by contradiction, assume that \( t = \operatorname{qid}_R M > n \). From the previous statement, we have \( \operatorname{Ext}_R^{t}(R/(\boldsymbol{x}), Z_s) = 0 \). Note that, by the choice of \( I \), it follows that \( t = \operatorname{id}_R Z_s \). Since \( Z_s \) is \(\operatorname{Ext}\)-finite by Lemma \ref{lemma:similar}, it follows from Lemma \ref{lemma:enochs} that   $ t = \sup \{ i : \operatorname{Ext}_R^i(k, Z_s) \neq 0 \}$. Consequently, we have \(\operatorname{Ext}_R^t(k, Z_s) \neq 0\). However, since \( \mathfrak{m} \in \operatorname{Ass}(R/(\boldsymbol{x})) \), there exists an inclusion \( 0 \rightarrow k \rightarrow R/(\boldsymbol{x}) \). This induces an exact sequence:
\[
\operatorname{Ext}^t_R(R/(\boldsymbol{x}), Z_s) \rightarrow \operatorname{Ext}_R^t(k, Z_s) \rightarrow 0.
\]

Consequently, \( \operatorname{Ext}_R^t(R/(\boldsymbol{x}), Z_s) = 0 \) implies \( \operatorname{Ext}_R^t(k, Z_s) = 0 \), which leads to a contradiction. Then, we conclude that $\operatorname{qid}_R M \leq n$.
\end{proof}
The following corollary is analogous to the corresponding result for Gorenstein-injective dimension \cite[Corollary 4.4]{ENOCHS}, as well as to the classical result for injective dimension \cite[Proposition 3.1.14]{BH}.
\begin{corollary}\label{crl:bom}
Let \((R, \mathfrak{m},k)\) be a local ring and let \(M\) be a non-zero \(\operatorname{Ext}\)-finite \(R\)-module of finite quasi-injective dimension. Then, for any maximal \(R\)-sequence \(\boldsymbol{x} \in \mathfrak{m}\), the following holds.
\begin{enumerate}
    \item $\operatorname{qid}_R M =0$ if and only if $\operatorname{Ext}_R^{i>0}(R/(\boldsymbol{x}),M)=0.$
    \item If $\operatorname{qid}_R M \neq 0$, then         $\operatorname{qid}_R M = \sup \lbrace i:\operatorname{Ext}_R^i(R/(\boldsymbol{x}),M ) \neq 0\rbrace$. 
\end{enumerate}
\end{corollary}
\begin{proof}
(1) Follows by (4) $\Leftrightarrow$ (1) of Theorem \ref{prop1}.

(2) Since $\operatorname{Ext}_R^i(R/(\boldsymbol{x}),M)=0$ for all $i> \operatorname{qid}_R M$ (see Theorem \ref{prop1} or Fact \ref{newfact}), we have that $t=\operatorname{qid}_R M \geq  \sup \lbrace i:\operatorname{Ext}_R^i(R/(\boldsymbol{x}),M ) \neq 0\rbrace$. By contradiction, assume that  $t>\sup \lbrace i:\operatorname{Ext}_R^i(R/(\boldsymbol{x}),M ) \neq 0\rbrace$. Hence, we have that  $\operatorname{Ext}_R^{i>t-1}(R/(\boldsymbol{x}),M)=0$. Therefore, by Theorem \ref{prop1} we have $\operatorname{qid}_R M \leq t-1$, which is a contradiction. Thus, $t=\sup \lbrace i:\operatorname{Ext}_R^i(R/(\boldsymbol{x}),M ) \neq 0\rbrace$.
\end{proof}

\begin{corollary}\label{ineq}
Let $R$ be a local ring and let $M$ be a non-zero $\operatorname{Ext}$-finite $R$-module of finite quasi-injective dimension. Then $\operatorname{qid}_R M \leq \operatorname{depth} R$.
\end{corollary}
\begin{proof}
Let $\boldsymbol{x} \in \mathfrak{m}$ be a maximal $R$-sequence. Since $\operatorname{pd}_RR/(\boldsymbol{x})=\operatorname{depth} R$, then the assertion follows by $(4) \Rightarrow (1)$ of Theorem \ref{prop1}.
\end{proof}
\begin{remark}
Let $R$ be a local ring. If $M$ is a finitely generated $R$-module with finite quasi-injective dimension, we can recover the Bass formula for quasi-injective dimension (\cite[Theorem 3.2]{Gheibi}) from Corollary \ref{crl:bom}. Indeed, since $M$ is finitely generated and $\operatorname{pd}_R R/(\boldsymbol{x}) = \operatorname{depth} R$, it follows from Corollary \ref{crl:bom} and \cite[p. 154, Lemma 1(iii)]{Matsu} that  
\begin{align*}
\operatorname{qid}_R M = \sup \lbrace i : \operatorname{Ext}_R^i(R/(\boldsymbol{x}), M) \neq 0 \rbrace = \operatorname{depth} R.
\end{align*}
\end{remark}

The following proposition completes the investigation into the quasi-injective dimension of quotients by regular elements and cosyzygy modules, which was initiated in Propositions 2.6(2) and 2.4(3) of \cite{Gheibi}, respectively. This result is established under the assumption that \( M \) is \(\operatorname{Ext}\)-finite.
\begin{proposition}\label{prop:formulas}
Let $(R,\mathfrak{m},k)$ be a local ring and let $M$ be a non-zero $\operatorname{Ext}$-finite $R$-module of finite quasi-injective dimension. Then

\begin{enumerate}
    \item If $y \in \mathfrak{m}$ is an element which is $R$-regular and $M$-regular, then $\operatorname{qid}_{R/(y)} M/yM=\sup \lbrace \operatorname{qid}_R M-1,0\rbrace$.
    \item  One has $\operatorname{qid}_R(\Omega_{-n}^RM) = \sup \lbrace \operatorname{qid}_R M-n,0 \rbrace $ for all $n \geq 0$.
\end{enumerate} 
\end{proposition}
\begin{proof}
(1) It is known that $\operatorname{qid}_{R/(y)} M/yM \leq \operatorname{qid}_R M$, as established in \cite[Proposition 2.6(2)]{Gheibi}. Consequently, $\operatorname{qid}_{R/(y)} M/yM<\infty$. Set $\overline{R}=R/(y)$ and $\overline{M}=M/yM$. For any $\overline{R}$-module $N$, we have the isomorphism $$\operatorname{Ext}_R^{i+1}(N,M) \cong \operatorname{Ext}_{\overline{R}}^i(N,\overline{M})$$ for all $i>0$, by \cite[p. 140, Lemma 2(i)]{Matsu}. Since $M$ is an $\operatorname{Ext}$-finite $R$-module, it follows that $\overline{M}$ is also an $\operatorname{Ext}$-finite $\overline{R}$-module. If $\operatorname{depth} R =1$, then $\operatorname{depth} R/(y)=0$ and $\operatorname{qid}_{R/(y)} M/yM=0$, by Corollary \ref{ineq}. Otherwise, consider $\boldsymbol{x}= y,x_1,\dots,x_n$ a maximal $R$-sequence in $\mathfrak{m}$ and $\overline{\boldsymbol{x}}=\overline{x_1},\dots,\overline{x_n} \in \overline{R}$. Note that $R/(\boldsymbol{x})$ is an $\overline{R}$-module and $$\operatorname{Ext}_R ^{i+1}(R/(\boldsymbol{x}),M) \cong \operatorname{Ext}_{\overline{R}}^{i}(R/(\boldsymbol{x}),\overline{M}) \cong \operatorname{Ext}^i_{\overline{R}}(\overline{R}/(\overline{\boldsymbol{x}}),\overline{M})$$ for all $i>0$. Therefore, Corollary \ref{crl:bom} yields the desired equality.

(2) The case \( n = 0 \) is trivial. By \cite[Proposition 2.4(3)]{Gheibi}, we have \(\operatorname{qid}_R(\Omega_{-1}^R M) < \infty\). Using \cite[Proposition 2.4(3)]{Gheibi}, the exact sequence 
$$ 
0 \rightarrow \Omega_{-n} M \rightarrow E^n(M) \rightarrow \Omega_{-(n+1)} M \rightarrow 0,$$ 
\noindent and proceeding by induction, it follows that \(\operatorname{qid}_R(\Omega_{-n}^R M) < \infty\) for all \( n \geq 0 \). 

Now, let $n>0$. For any \( R \)-module \( N \), we have the isomorphism  
\[
\operatorname{Ext}_R^i(N, \Omega_{-n}^R M) \cong \operatorname{Ext}_R^{i+n}(N, M)  
\]
for all \( i > 0 \). Since \( M \) is \(\operatorname{Ext}\)-finite, it follows that \(\Omega_{-n}^R M\) is also \(\operatorname{Ext}\)-finite. In particular, if $\boldsymbol{x} \in \mathfrak{m}$ is a maximal $R$-regular sequence, then $$\operatorname{Ext}_R^i(R/(\boldsymbol{x}), \Omega_{-n}^R M) \cong \operatorname{Ext}_R^{i+n}(R/(\boldsymbol{x}), M).$$ Hence, by applying the homological characterization of \(\operatorname{qid}_R M\) given in Corollary \ref{crl:bom}, we obtain the desired equality.
\end{proof}

The following proposition improves \cite[Corollary 3.3]{Gheibi}, since we are not considering that $R$ is local or that $M$ is finitely generated.

\begin{proposition}\label{prop:comparativo}
Let $M$ be an $R$-module of finite quasi-injective dimension. Then $\operatorname{qid}_R M \leq \operatorname{id}_R M$ and equality holds when $\operatorname{id}_R M$ is finite.
\end{proposition}
\begin{proof}
Let $\operatorname{id}_RM<\infty$. If $\operatorname{qid}_RM>0$, then using the Chouinard's formula for quasi-injective dimension (see \cite[Theorem 3.3]{Tri}) and the Classical Chouinard formula (\cite[Corollary 3.1]{Chouinard}), we have
\begin{align*}
\operatorname{qid}_R M = \sup \lbrace \operatorname{depth} R_{\mathfrak{p}}- \operatorname{width}_{R_{\mathfrak{p}}} M_{\mathfrak{p}} \mid \mathfrak{p} \in \operatorname{Spec}R \rbrace=\operatorname{id}_RM.
\end{align*}
Now, assume that $\operatorname{qid}_RM=0$. Since $\operatorname{id}_RM< \infty$, then for each $\mathfrak{p} \in \operatorname{Spec} R$, we have $\operatorname{Ext}_R^{i \gg 0}(R/\mathfrak{p},M)=0$. Since $\operatorname{qid}_R M =0$, then $\operatorname{Ext}_R^1(R/\mathfrak{p},M)=0$ for all $\mathfrak{p} \in \operatorname{Spec} R$, by Fact \ref{newfact}. Thus $\operatorname{id}_R M =0$, by \cite[Corollary 3.1.12]{BH}.
\end{proof}

Recall that for \(x \in R\), an \(R\)-module \(M\) is called \(x\)-\textit{divisible} if \(xM = M\).

\begin{lemma}\label{lemma:divisible}
Let $M$ be an $R$-module with $\operatorname{qid}_R M = 0$. Then
\begin{enumerate}
    \item If $x \in R$ is $R$-regular, then $M$ is $x$-divisible. 
    \item If $M$ is not injective, then $\operatorname{id}_R M=\infty$. 
\end{enumerate}
\end{lemma}
\begin{proof}
(1) Let $\overline{R}=R/(x)$. Consider the exact sequence $0 \rightarrow R \xrightarrow{.x} R \rightarrow \overline{R} \rightarrow 0$ and applying the functor $\operatorname{Hom}_R (-,M)$ we obtain the exact sequence $$0 \rightarrow \operatorname{Hom}_R (\overline{R},M) \rightarrow M \xrightarrow{.x} M \rightarrow0$$
using Fact \ref{newfact} to obtain that $\operatorname{Ext}_R^1(\overline{R},M)=0$,  since  $\operatorname{qid}_R M =0 $ and $\operatorname{pd}_R \overline{R} \leq 1$. Thus, $M=xM$.

(2) Follows by Proposition \ref{prop:comparativo}.
\end{proof}
\begin{remark}
If $R$ is local and $\operatorname{depth} R> 0$, then every non-zero $R$-module with $\operatorname{qid}_R M=0$ is not finitely generated, by Lemma \ref{lemma:divisible}(1).
\end{remark}
\begin{lemma}\label{isomorphisms}
(\cite[Lemma 2.1]{ENOCHS}) Let $x$ be an $R$-regular element and let $\overline{R}=R/(x)$. Let $M$ be a $x$-divisible $R$-module, and $N$ be an $\overline{R}$-module. Then
\begin{enumerate}
    \item $\operatorname{Ext}_R^i(N,M)   \cong  \operatorname{Ext}_{\overline{R}}^i(N,\operatorname{Hom}_{R}(\overline{R},M)).$
    \item $\operatorname{Ext}_R^{i+1}(M,N) \cong \operatorname{Ext}_{\overline{R}}^i(\operatorname{Hom}_R(\overline{R},M),N).$
    \item $\operatorname{Tor}_{i+1}^R(M,N)\cong \operatorname{Tor}_i^{\overline{R}} (\operatorname{Hom}_R(\overline{R},M),N).$
\end{enumerate}
\end{lemma}
\begin{theorem}\label{lemma:provar}
Let $(R,\mathfrak{m},k)$ be a local ring and let $M$ be a non-injective $\operatorname{Ext}$-finite $R$-module with $\operatorname{qid}_R M=0$. If $x \in \mathfrak{m}$ is $R$-regular and $\overline{R}=R/(x)$, then $\operatorname{Hom}_{R}(\overline{R},M)$ is a non-injective $\operatorname{Ext}$-finite $\overline{R}$-module and $\operatorname{qid}_{\overline{R}} (\operatorname{Hom}_{R}(\overline{R},M))=0$. 
\end{theorem}
\begin{proof}
By Lemma \ref{lemma:divisible}(1), we have that $M$ is $x$-divisible. Therefore, for every finitely generated $\overline{R}$-module $N$, which is also a finitely generated $R$-module, Lemma \ref{isomorphisms}(1) establishes the isomorphism $\operatorname{Ext}_R^i(N,M) \cong \operatorname{Ext}_{\overline{R}}^i(N,\operatorname{Hom}_{R}(\overline{R},M))$ for all $i \geq 1$. Thus, since $M$ is $\operatorname{Ext}$-finite, the same holds for $\operatorname{Hom}_R(\overline{R},M)$.

Now, we prove that $\operatorname{qid}_{\overline{R}} (\operatorname{Hom}_R(\overline{R},M)) \leq \operatorname{qid}_R M=0$. If $\operatorname{Hom}_{R}(\overline{R},M)=0,$ then $\operatorname{qid}_{\overline{R}} (\operatorname{Hom}_{R}(\overline{R},M))=-\infty$ and this inequality is trivial. So, in order to prove this inequality, we may assume that $\operatorname{Hom}_R(\overline{R},M) \neq 0$. Let $I$ be a quasi-injective resolution of $M$. That is, $H_q(I)\cong M^{\oplus b_q}$ for some $b_q \geq 0$, for each $q \in \mathbb{Z}$. We claim that $\operatorname{Hom}_R(\overline{R},I)$ is a quasi-injective resolution of the $\overline{R}$-module $\operatorname{Hom}_R(\overline{R},M)$. Indeed, by \cite[Lemma 2.5]{Gheibi}, there exists a convergent spectral sequence 
\begin{align*}
   E_2^{p,q} = \operatorname{Ext}_R^p(\overline{R},H_q(I)) \Longrightarrow H_{q-p}(\operatorname{Hom}_R(\overline{R},I)).
\end{align*}
By \cite[Lemma 3.1.6]{BH}, \(\operatorname{Hom}_R(\overline{R}, I)\) forms a complex of injective \(\overline{R}\)-modules. Since \(x\) is \(R\)-regular and $\operatorname{qid}_R M=0$, using Fact \ref{newfact} we have $$\operatorname{Ext}_R^{i}(\overline{R},H_q(I))\cong \operatorname{Ext}_R^i(\overline{R},M)^{\oplus b_q} =0,$$ for all $i>0$ for each $q \in \mathbb{Z}$. Then
\begin{align*}
E^{p,q}_{\infty} \cong \begin{cases}
\operatorname{Hom}_R(\overline{R},H_q(I)) & \text{if } p=0 \\
0 & \text{otherwise.}
\end{cases}
\end{align*}
Therefore, we obtain 
\begin{align*}
 \operatorname{Hom}_R (\overline{R}, M) ^{\oplus b_q} 
  \cong  \operatorname{Hom}_R( \overline{R}, H_q(I)) \cong H_q(\operatorname{Hom}_R(\overline{R},I)).
\end{align*}

That is, $\operatorname{Hom}_R(\overline{R},I)$ is a quasi-injective resolution of the $\overline{R}$-module $\operatorname{Hom}_R(\overline{R},M)$ and $\operatorname{hinf} (\operatorname{Hom}_R(\overline{R},I))=\operatorname{hinf} I$ . Thus we obtain $\operatorname{qid}_{\overline{R}} (\operatorname{Hom}_R(\overline{R},M)) \leq \operatorname{qid}_R M=0$, as asserted. 

To complete this proof, we need to show that $\operatorname{Hom}_R(\overline{R},M) \neq 0$ and that $\operatorname{Hom}_R (\overline{R},M)$ is not injective.

Since $M$ is not injective, then $\operatorname{id}_R M = \infty$, by Lemma \ref{lemma:divisible}(2). Thus, the $i$-th Bass number of $M$ is non-zero for all $i \geq \dim R$ by \cite[Proposition 2.5]{ENOCHS} and it means that $E(k)$ is a summand of $E^i(M)$. Then $\operatorname{Hom}_R(\overline{R},E(k))$ is a summand of $\operatorname{Hom}_R (\overline{R},E^i(M))$. Following the method presented in the proof of \cite[Lemma 2.1]{ENOCHS}, we observe that if $0 \rightarrow M \rightarrow E(M) \rightarrow E^1(M) \rightarrow \cdots$ is the minimal injective resolution of $M$, then $$0 \rightarrow\operatorname{Hom}_R(\overline{R},M) \rightarrow \operatorname{Hom}_R (\overline{R},E(M)) \rightarrow \operatorname{Hom}_R(\overline{R},E^1(M)) \rightarrow \cdots$$ is the minimal $\overline{R}$-injective resolution of $\operatorname{Hom}_R(\overline{R},M)$, since $M$ is $x$-divisible and thus $\operatorname{Ext}_R^{i >0}(\overline{R},M)=0$. Therefore, $E^i_{\overline{R}}(\operatorname{Hom}_R (\overline{R},M)) = \operatorname{Hom}_R (\overline{R},E^i(M)).$ Given that $\operatorname{Hom}_R (\overline{R},E(k))\neq 0$, it follows that $E^i_{\overline{R}}(\operatorname{Hom}_R(\overline{R},M)) \neq 0$ for all $i \geq \dim R$. Thus $\operatorname{Hom}_R (\overline{R},M) \neq 0$ and $\operatorname{id}_{\overline{R}} (\operatorname{Hom}_R (\overline{R},M))=\infty$, as desired. 
\end{proof}

\section{Quasi-injective dimension and 
width of a module}
In this section, we are now able to prove our main theorem. This theorem can be compared with previously established results concerning Gorenstein-injective dimension and injective dimension, such as \cite[Theorem 4.8]{ENOCHS}, \cite[Theorem 2.5]{SEZEDEH}, \cite[Theorem 2.4]{SHARIF}, and \cite[Theorem 1.6]{YASSEMI}. Let $(R,\mathfrak{m},k)$ be a local ring and let $M$ be an $R$-module, we recall that $\operatorname{width}_R M = \inf \lbrace i : \operatorname{Tor}_i^R(k,M) \neq 0 \rbrace$. 

\begin{theorem}\label{theorem1}
Let $(R,\mathfrak{m},k)$ be a Cohen-Macaulay local ring and let $M$ be an $\operatorname{Ext}$-finite $R$-module with $\operatorname{width}_R M<\infty$ (e.g., if $M$ is $\operatorname{Ext}$-finite and non-injective). If $M$ is of finite quasi-injective dimension, then $$\operatorname{qid}_R M + \operatorname{width}_R M = \operatorname{depth} R.$$ 
\end{theorem}

\begin{proof}
If $M$ is injective, then $0=\operatorname{id}_R M \operatorname=\operatorname{qid}_RM$ and the equality follows by \cite[Theorem 1.6]{YASSEMI}. Thus, we may assume the additional hypothesis that \( M \) is non-injective. The proof of this theorem proceeds by induction on \( \operatorname{qid}_R M \).

Let $\operatorname{qid}_RM=0$. In this case, we proceed by induction on $d=\operatorname{depth} R$. If $d=0$, then the maximal ideal $\mathfrak{m}$ is nilpotent. Then, $R/\mathfrak{m} \otimes_ R M \cong M/\mathfrak{m}M \neq 0$ and so $\operatorname{width}_R M =0$. Thus $\operatorname{width}_R M = \operatorname{depth} R$, as desired. 
Assume that $d \geq 1$ and let $x \in \mathfrak{m}$ be an $R$-regular element. Set $\overline{R}=R/(x)$. Then  $\operatorname{Hom}_{R}(\overline{R},M)$ is a non-injective $\operatorname{Ext}$-finite $\overline{R}$-module and $\operatorname{qid}_{\overline{R}} (\operatorname{Hom}_{R}(\overline{R},M))=0$, by Theorem \ref{lemma:provar}. Thus, by induction, we have that $\operatorname{width}_{\overline{R}} ( \operatorname{Hom}_R(\overline{R}, M)) = d-1$. Note that $M$ is $x$-divisible, by Lemma \ref{lemma:divisible}(1). Therefore, the isomorphism given in Lemma $\ref{isomorphisms}(3)$ provides the equality $\operatorname{width}_{\overline{R}} (\operatorname{Hom}_R(\overline{R},M))=\operatorname{width}_R M -1$. Thus $\operatorname{width}_R M =d$.

Now, assume that $\operatorname{qid}_R M=1$ and then $ d\geq 1$, by Corollary \ref{ineq}. In this case, $\Omega_{-1}^R M= E(M)/M$ satisfies $\operatorname{qid}_R ( \Omega_{-1}^R M) =0$, by Proposition \ref{prop:formulas}(2). If $\Omega_{-1}^R M$ is injective, then $\operatorname{qid}_RM=\operatorname{id}_R M=1$ and the formula follows by \cite[Theorem 1.6]{YASSEMI}. Thus, we may assume that $\Omega_{-1}^RM$ is non-injective and then $\operatorname{width}_R (\Omega_{-1}^R M) < \infty$, since it is $\operatorname{Ext}$-finite (see Theorem \ref{finite}). Therefore $\operatorname{width}_R (\Omega_{-1}^R M)=\operatorname{depth} R$, by the previous case. Theorem \ref{finite} provides the following inequality: \(\operatorname{width}_R M \leq d\). If \(\operatorname{width}_R M \leq d - 1\), then
\[
    \operatorname{width}_R M + 1 = \operatorname{width}_R (\Omega_{-1}^R M) = d,
\]
by \cite[Lemma 4.6]{ENOCHS}, and the proof of this case is complete. Therefore, by contradiction, assume that \(\operatorname{width}_R M = d\). Consider the notation introduced in Lemma \ref{lemma:similar} and set $s=\operatorname{hinf} I$. By the choice of $I$, we have that $Z_s$ is $\operatorname{Ext}$-finite (see Lemma \ref{lemma:similar}), $\operatorname{id}_R Z_s=1$ and then $\operatorname{width}_R Z_s< \infty$, by Theorem \ref{finite}.  Therefore, by \cite[Theorem 1.6]{YASSEMI}, we have $\operatorname{width}_R Z_s = d -1$. There are exact sequences
\begin{align}\label{seq3}
\begin{cases}
0 \rightarrow Z_i \rightarrow I_i \rightarrow B_{i-1} \rightarrow 0 \\  0 \rightarrow B_i \rightarrow Z_i \rightarrow H_i(I) \rightarrow 0
\end{cases}  
\quad (i \in \mathbb{Z})
\end{align}
with $H_i(I) \cong M^{\oplus b_i}$, for some $b_i \geq 0$. Since  we are assuming $\operatorname{width}_RM=d$, then $\operatorname{Tor}_i^R(k,M)=0$ for $i=0,1,\dots,d-1$ and the sequence $0 \rightarrow B_s \rightarrow Z_s \rightarrow M^{\oplus b_s} \rightarrow0 $ induces the following exact sequence: 
\begin{align*}
\cdots \rightarrow \operatorname{Tor}_d^R 
 (k,M)^{\oplus b_s} \rightarrow \operatorname{Tor}_{d-1}^R(k,B_s) \rightarrow \operatorname{Tor}_{d-1}^R (k,Z_s) \rightarrow0 
\end{align*}
and then $\operatorname{Tor}_{d-1}^R(k,B_s) \neq 0$, as $\operatorname{Tor}_{d-1}^R(k,Z_s) \neq 0$ since we have $\operatorname{width}_R Z_s = d -1$. Considering the sequence $0 \rightarrow Z_{s+1} \rightarrow I_{s+1} \rightarrow B_s \rightarrow 0$, one can see that $\operatorname{width}_R I_{s+1} $ is $d$ or $\infty$ (see Proposition \ref{vanishing}). In either case, $\operatorname{Tor}_i^R(k,I_{s+1})=0$ for $i=0,1,\dots,d-1$. Thus, this sequence induces the following exact sequence:
\begin{align*}
    0 \rightarrow\operatorname{Tor}_{d-1}^R (k,B_s) \rightarrow \operatorname{Tor}_{d-2}^R(k,Z_{s+1})
\end{align*}
and then $\operatorname{Tor}_{d-2}^R(k,Z_{s+1}) \neq 0$, since $\operatorname{Tor}_{d-1}^R(k,B_s) \neq 0$. Using the exact sequences (\ref{seq3}) and repeating this argument a finite number of times, we have $k \otimes_R B_{s+d-1} \neq 0$. Also, applying Proposition \ref{vanishing} once more, we see that $\operatorname{width}_R I_{s+d}$ is $d$ or $\infty$ and then $k \otimes_R I_{s+d} =0$. Then, by tensoring the exact sequence $$0 \rightarrow Z_{s+d} \rightarrow I_{s+d} \rightarrow B_{s+d-1} \rightarrow 0$$ by $k$ we have a contradiction, since $k \otimes_R  B_{s+d-1} \neq 0$. Therefore $\operatorname{width}_R M \leq d-1$.

Finally, assume that $\operatorname{qid}_R M >1$. Then $\operatorname{qid}_R (\Omega_{-1}^RM) =\operatorname{qid}_R M-1 \geq 1$, by Proposition \ref{prop:formulas}(2). Note that $\Omega_{-1}^R M$ is Ext-finite and non-injective. Then $\operatorname{width}_R (\Omega_{-1}^R M) < \infty$, by Theorem \ref{finite}. Thus, by induction, we have 
\begin{align*}
    \operatorname{qid}_R M-1 = \operatorname{qid}_R (\Omega_{-1}^RM) = d- \operatorname{width}_R( \Omega_{-1}^R M).
\end{align*}
Since $\operatorname{width}_R (\Omega_{-1}^R M) < d$ and $\operatorname{Tor}_i^R(k,E(M))=0$ for $i=0,1,\dots,d-1$ (see Proposition \ref{vanishing}), one can see that $\operatorname{width}_R (\Omega_{-1}^R M) = \operatorname{width}_R M +1$, and thus $\operatorname{qid}_R M + \operatorname{width}_ RM = d$.
\end{proof}

\begin{example}
Let $(R,\mathfrak{m},k)$ be a non-regular Cohen-Macaulay local ring. Therefore $\operatorname{id}_R k =\infty$ (see \cite[Exercise 3.1.26]{BH}) and $\operatorname{qid}_R k = \operatorname{depth}R$ (See \cite[Proposition 2.8(1) and Theorem 3.2]{Gheibi}). Each cosyzygy module $\Omega_{-n}^R k$ is $\operatorname{Ext}$-finite and non-injective and then $\operatorname{width}_R (\Omega_{-n}^Rk) <\infty$, by Theorem \ref{finite}. Using Proposition \ref{prop:formulas}(2) and Theorem \ref{theorem1}, we have the following formulas: 
\begin{align*}
\operatorname{qid}_R (\Omega_{-n}^Rk) & = \begin{cases}
 \operatorname{depth} R- n  ,  & \text{ if } n\leq\operatorname{depth } R \\
0, & \text{ if } n > \operatorname{depth} R
 \end{cases}
 \end{align*}
 and
 \begin{align*}
\operatorname{width}_R (\Omega_{-n}^Rk) & =\begin{cases}
 n ,  & \text{ if } n\leq\operatorname{depth } R \\
 \operatorname{depth} R, & \text{ if } n > \operatorname{depth} R.
\end{cases}
\end{align*}
This example shows that the inequality in Corollary \ref{ineq} can be strict when $M$ is not finitely generated, even when $\operatorname{id}_R M = \infty$. 
\end{example}

\begin{agra}
The authors would like to thank Victor D. Mendoza Rubio for his valuable comments during the preparation of this manuscript. The authors would also like to thank the referee for the careful reading and for his/her comments and suggestions, which have improved the quality of the paper. 
\end{agra}

\noindent
\textbf{Disclosure statement.} No potential conflict of interest was reported by the authors. 

\begin{fund}
The first author was supported by  S\~ao Paulo Research Foundation (FAPESP) under grant 2019/21181-0. The second author was supported by  S\~ao Paulo Research Foundation (FAPESP) under grant 2022/12114-0.   
\end{fund}

\end{document}